\documentclass[12pt]{article}
\usepackage{tikz}
\usetikzlibrary{calc,through,backgrounds}
\usepackage{multicol}

\usepackage{bm,amssymb,amsthm,amsfonts,amsmath,latexsym,cite,color, psfrag,graphicx,ifpdf,pgfkeys,pgfopts,xcolor,extarrows}

\newtheorem{theo}{Theorem}[section]

\newtheorem{coro}[theo]{Corollary}

\allowdisplaybreaks

\makeatletter \@addtoreset{equation}{section}
\@addtoreset{theo}{}\makeatother

\usepackage{hyperref}

\def\qed{\hfill \rule{4pt}{7pt}}
\def\pf{\noindent {\it Proof.} }
\def\S{  \mathfrak{S}}
\def\N{  \mathrm{N}}
\def\Min{  \mathrm{Min}}
\def\Max{  \mathrm{Max}}
\def\GL{   \mathrm{GL}  }
\def\Char{   \mathrm{char}  }
\def\tr{   \mathrm{tr}  }

\begin{document}
\begin{center}
{\Large\bf Upper Bounds  of Schubert Polynomials}

\vskip 6mm
{\small  Neil J.Y. Fan$^1$ and Peter L. Guo$^2$ }

\vskip 4mm
$^1$Department of Mathematics\\
Sichuan University, Chengdu, Sichuan 610064, P.R. China
\\[3mm]

$^{2}$Center for Combinatorics, LPMC\\
Nankai University,
Tianjin 300071,
P.R. China

\vskip 4mm

$^1$fan@scu.edu.cn, $^2$lguo@nankai.edu.cn
\end{center}

\begin{abstract}

Let $w$ be a permutation  of   $\{1,2,\ldots,n \}$, and let $D(w)$ be the Rothe diagram of $w$.
The Schubert polynomial $\S_w(x)$ can be realized as the dual
character of the flagged Weyl module associated to  $D(w)$.
This implies  a coefficient-wise inequality
\[\Min_w(x)\leq \S_w(x)\leq \Max_w(x),\]
where both $\Min_w(x)$ and $\Max_w(x)$ are  polynomials
 determined by   $D(w)$.
Fink,  M\'esz\'aros and  St.\,Dizier found that $\S_w(x)$ equals the lower bound $\Min_w(x)$ if and only if $w$ avoids
twelve permutation patterns. In this paper, we show that
$\S_w(x)$ reaches the upper bound $\Max_w(x)$ if and only if $w$ avoids  two permutation  patterns 1432 and 1423.
Similarly, for any given composition $\alpha\in \mathbb{Z}_{\geq 0}^n$, one can define  a lower bound $\Min_\alpha(x)$ and an
upper bound  $\Max_\alpha(x)$ for  the key polynomial $\kappa_\alpha(x)$.
Hodges and Yong established  that   $\kappa_{\alpha}(x)$ equals $\Min_\alpha(x)$  if and only if $\alpha$ avoids five composition patterns. We show that $\kappa_{\alpha}(x)$ equals $\Max_\alpha(x)$
if and only if $\alpha$ avoids a single composition pattern $(0,2)$.
As an application, we obtain that when $\alpha$ avoids  $(0,2)$, the key polynomial
$\kappa_{\alpha}(x)$  is  Lorentzian, partially  verifying a conjecture of
Huh,   Matherne,   M\'esz\'aros and  St.\,Dizier.

\end{abstract}

\section{Introduction}

The Schubert polynomials  $\S_w(x)=\S_w(x_1,\ldots,x_n)$  indexed by  permutations
 $w$ of $[n]=\{1,2,\ldots,n \}$ were introduced by Lascoux and Sch\"utzenberger \cite{Las},
 representing cohomology classes of Schubert cycles in flag varieties.
For combinatorial constructions of Schubert polynomials, see
for example \cite{Ber,Bil,Knu,Lam,Win}.

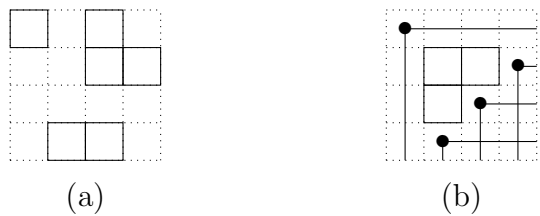
\begin{figure}[h]
\begin{center}
\begin{tikzpicture}

\def\rectanglepath{-- +(5mm,0mm) -- +(5mm,5mm) -- +(0mm,5mm) -- cycle}

\draw [step=5mm,dotted] (0mm,0mm) grid (20mm,20mm);
\draw (0mm,15mm) \rectanglepath;\draw (10mm,15mm) \rectanglepath;
\draw (10mm,10mm) \rectanglepath;\draw (15mm,10mm) \rectanglepath;
\draw (5mm,0mm) \rectanglepath;\draw (10mm,0mm) \rectanglepath;

\draw[dotted](50mm,0mm)--(50mm,20mm);
\draw [step=5mm,dotted] (50mm,0mm) grid (70mm,20mm);
\draw (55mm,10mm) \rectanglepath;\draw (60mm,10mm) \rectanglepath;
\draw (55mm,5mm) \rectanglepath;
\node at (52.5mm,17.5mm) {$\bullet$};
\node at (67.5mm,12.5mm) {$\bullet$};
\node at (62.5mm,7.5mm) {$\bullet$};
\node at (57.5mm,2.5mm) {$\bullet$};
\draw(52.5mm,17.5mm)--(70mm,17.5mm);\draw(52.5mm,17.5mm)--(52.5mm,0mm);
\draw(67.5mm,12.5mm)--(70mm,12.5mm);\draw(67.5mm,12.5mm)--(67.5mm,0mm);
\draw(62.5mm,7.5mm)--(70mm,7.5mm);\draw(62.5mm,7.5mm)--(62.5mm,0mm);
\draw(57.5mm,2.5mm)--(70mm,2.5mm);\draw(57.5mm,2.5mm)--(57.5mm,0mm);

\node at (10mm,-5mm) {(a)};\node at (60mm,-5mm) {(b)};

\end{tikzpicture}
\end{center}
\vspace{-6mm}
\caption{(a) is a diagram,   (b) is the Rothe diagram of $w=1432$.}
\label{diagram}
\end{figure}

Kra\'skiewicz and  Pragacz \cite{Kra,Kra-2} proved that $\S_w(x)$ equals the
dual character of the flagged Weyl module associated to the Rothe diagram $D(w)$ of $w$.
Given a  permutation $w=w_1w_2\cdots w_n$ of $[n]$, the Rothe diagram
$D(w)=(D(w)_1,D(w)_2,\ldots, D(w)_n)$
of $w$ is  defined by
\[D(w)_j=\{i\colon w_i>j,\ i<w^{-1}_j\},\ \ \ \text{ where $1\leq j\leq n$.}\]
In general, a diagram means an ordered list $D=(D_1,D_2,\ldots, D_n)$ of $n$ subsets of $[n]$.
A diagram $D$ can   be viewed as a collection of boxes of
an $n\times n$ grid, this is, $D_j$ consists of the boxes $(i,j)$ in row $i$ and column $j$ where $i\in D_j$.
Here  the row indices increase from top to bottom, and the column indices increase  from left to right.
  For example, Figure \ref{diagram}(a) represents
   the diagram $(\{1\}, \{4\}, \{1,2,4\},\{2\})$.
When viewed as a subset of an $n\times n$ grid, the Rothe diagram $D(w)$
can be obtained   by removing the boxes that are to the right of  $(i, w_i)$
or below $(i, w_i)$. Figure \ref{diagram}(b) is
the Rothe diagram of $1432$.

\vspace{2mm}

For two diagrams $C=(C_1,\ldots,C_n)$ and $D=(D_1,\ldots,D_n)$, write $C\leq D$ if
$C_j\leq D_j$ for every $1\leq j\leq n$, where $C_j\leq D_j$ means that
\begin{itemize}
\item[(1)] $|C_j|=|D_j|$;

\item[(2)] for $1\leq k\leq |C_j|$,
the $k$-th least element of $C_j$ is less than or equal to the $k$-th least element of $D_j$.
\end{itemize}
It is worth mentioning that the set $\{C_j\colon C_j\leq D_j\}$ forms the basis of the
Schubert matroid corresponding to $D_j$, see for example \cite{Ard,Bon,Fan}.

Write $x^D$ for the monomial generated by a diagram $D$:
\[x^D=\prod_{ j=1}^n \prod_{i\in D_j}x_i.\]
Recall that $D(w)$ is the Rothe diagram of a permutation $w$. Denote
\[
\Min_w(x)=\sum_{x^\beta\in \{x^C\colon C\leq D(w)\}}x^\beta
\ \ \ \ \ \ \text{and}\ \ \ \ \ \  \Max_w(x)=\sum_{  C\leq D(w) }x^C.
\]
Note that   the coefficient of each monomial
appearing in $\Min_w(x)$  equals one.

Given a diagram $D$, one can construct the flagged Weyl module $\mathcal{M}_D$ of the group $B$ of
invertible upper-triangular $n\times n$ matrices over $\mathbb{C}$ \cite{Fin-1,Fin,Kra,Kra-2,Mag},
see  Section \ref{Weyl} for  detailed descriptions.
Kra\'skiewicz and  Pragacz \cite{Kra, Kra-2} showed that $\S_w(x)$ equals the dual character
of   $\mathcal{M}_{D(w)}$.
As a consequence, one has the following coefficient-wise inequality:
\begin{align}\label{FG-1}
\Min_w(x)\leq \S_w(x)\leq \Max_w(x),
\end{align}
where, for two polynomials
$f(x)=\sum_{\beta}a_\beta\, x^\beta$ and $g(x)=\sum_{\beta}b_\beta\, x^\beta$ in
$\mathbb{Z}[x_1,\ldots,x_n]$,
$f(x)\leq g(x)$ means that $a_\beta\leq b_\beta$ for any $\beta$.
The above inequality \eqref{FG-1} will also be explained in Section \ref{Weyl}.

Fink,  M\'esz\'aros and  St.\,Dizier \cite[Theorem 1.1]{Fin} proved  that $\S_w(x)$ attains
 the lower bound $\Min_w(x)$   if and only if $w$ avoids
 twelve permutation patterns:
 12543, 13254,
13524, 13542, 21543, 125364, 125634, 215364, 215634, 315264, 315624,  315642.
For a permutation $w$ of $[n]$ and a permutation $\pi$ of $[m]$ with $n\geq m$, we say that $w$ avoids the pattern $\pi$
if there do not exist subsequences in $w$ of length $m$ that are order isomorphic to $\pi$.
A Schubert polynomial  that reaches the lower bound is  called  a zero-one Schubert polynomial
in \cite{Fin},
which can also be generated exactly by the lattice points in its associated Newton polytope  \cite{Fin-1}.

Our first result  provides a characterization of when  $\S_w(x)$
reaches  the upper bound.

\begin{theo}\label{main}
The Schubert polynomial $\S_w(x)$ equals $\Max_w(x)$ if and only if $w$ avoids the
patterns 1432 and 1423.
\end{theo}

Huh,   Matherne,   M\'esz\'aros and  St.\,Dizier
\cite[Conjecture 15]{Huh} conjectured that for any permutation $w$,
the normalized Schubert
polynomial  $\N(\S_w(x))$ is Lorentzian, where $\N$ is a linear operator defined by
\[\N(x^{\mu})=\frac{x^{\mu}}{\mu!}=\frac{x_1^{\mu_1}\cdots x_n^{\mu_n}}{\mu_1!\cdots\mu_n!},\quad
\text{ for $\mu=(\mu_1,\ldots,\mu_n)\in\mathbb{Z}_{\geq 0}^n$}.\]
We refer the reader  to  \cite{Bra} or \cite[Definition 5]{Huh}
 for several  equivalent definitions of Lorentzian polynomials.
It should be pointed out that, as an important consequence of the Lorentzian property,
the coefficients of a Lorentzian polynomial are log-concave.
Using Theorem \ref{main} combined with results in \cite{Bra,Fin-1},
Huh,   Matherne,   M\'esz\'aros and  St.\,Dizier
\cite[Proposition 17]{Huh} confirmed the above conjecture for permutations
avoiding 1432 and 1423. In fact, \cite[Proposition 17]{Huh} proved  that
when $w$ avoids 1432 and 1423, the Schubert polynomial $\S_w(x)$ is Lorentzian. This is  stronger   because
the Lorentzian property of a polynomial $f(x)$  implies that of  $\N(f(x))$  \cite[Corollary 3.7]{Bra}.
As noted below   \cite[Proposition 17]{Huh},
the Schubert polynomials $\S_{1432}(x)$ and $\S_{1423}(x)$ are not Lorentzian.

We also remark  that  Gire \cite{Gir} showed that
  the number of permutations of $[n]$ avoiding   2341 and 3241 is the large Schr\"{o}der number
  $r_{n-1}$  (see also Kremer \cite{Kremer}), which can be defined via
  the following generating function
  \[\sum_{r\geq 0}  r_nx^n=\frac{1-x-\sqrt{1-6x+x^2}}{2x}.\]
  The first few values of $r_n$ are $1, 2, 6, 22, 90, 394, 1806, 8558, 41586, \ldots$.
Reversing the order of permutations leads to a bijection between permutations of $[n]$
 avoiding  1432 and 1423 and    permutations of $[n]$
 avoiding  2341 and 3241. Thus the number of permutations of $[n]$ avoiding  1432 and 1423 is equal to
  $r_{n-1}$.

Using analogous arguments, we can characterize  when   key polynomials reach  their upper bounds.
Key polynomials $\kappa_\alpha(x)$
associated to compositions $\alpha\in \mathbb{Z}_{\geq 0}^n$ (also called Demazure characters)
 are characters of the Demazure
modules for the general linear groups \cite{Dem-1,Dem-2,Fin-1}.
Their combinatorial properties were initially investigated by Lascoux and Sch\"utzenberger
\cite{Las-1}.
It is known that every Schubert polynomial is a positive sum
 of key polynomials
  \cite{Las-2,Rei}.
It is also worth mentioning  that
    $\kappa_\alpha(x)$ can be
 realized as a   specialization of
the nonsymmetric Macdonald polynomial $E_\alpha(x; q, t)$
 at $q=t=0$  \cite{Ass,Ion}.

\begin{figure}[h]
\begin{center}
\begin{tikzpicture}[scale=0.5]
      \draw [-, shift={(0,0)}] (0,0)--(0,4);
      \filldraw [fill=gray, fill opacity=0, shift={(0,0)}] (0,0) rectangle (1,1);
      \filldraw [fill=gray, fill opacity=0, shift={(0,0)}] (1,0) rectangle (2,1);
      \filldraw [fill=gray, fill opacity=0, shift={(0,0)}] (0,2) rectangle (1,3);
      \filldraw [fill=gray, fill opacity=0, shift={(0,0)}] (1,2) rectangle (2,3);
      \filldraw [fill=gray, fill opacity=0, shift={(0,0)}] (2,2) rectangle (3,3);
      \filldraw [fill=gray, fill opacity=0, shift={(0,0)}] (0,3) rectangle (1,4);
\end{tikzpicture}
\end{center}
\vspace{-.3cm}
\caption{The skyline diagram of $(1,3,0,2)$. }\label{Fig-skyline}
\end{figure}
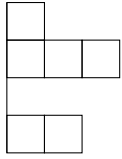

The key polynomial   $\kappa_\alpha(x)$  is equal to the dual character
of the  flagged Weyl module $\mathcal{M}_{D(\alpha)}$
associated to the skyline diagram $D(\alpha)$ of $\alpha$ \cite{Dem-1,Dem-2,Fin-1}.
Recall that the skyline diagram
  $D(\alpha)$
  consists of the first $\alpha_i$ boxes in row $i$.
 For example,  Figure \ref{Fig-skyline} depicts  the skyline diagram of  $(1,3,0,2)$.

Write
\[
\Min_\alpha(x)=\sum_{x^\beta\in \{x^C\colon C\leq D(\alpha)\}}x^\beta
\ \ \ \ \ \ \text{and}\ \ \ \ \ \  \Max_\alpha(x)=\sum_{  C\leq D(\alpha) }x^C.
\]
Similar to Schubert polynomials, one has the following inequality:
\begin{align}\label{FG-1x}
\Min_\alpha(x)\leq \kappa_\alpha(x)\leq \Max_\alpha(x).
\end{align}
Recently, Hodges and   Yong \cite{Hod} established a pattern avoidance characterization
of when $\kappa_\alpha(x)$ equals the lower bound $\Min_\alpha(x)$,
which they call multiplicity-free key polynomials.
For  two compositions $\alpha=(\alpha_1,\ldots, \alpha_n)$ and
$\beta=(\beta_1,\ldots, \beta_m)$ with $n\geq m$, we say that $\alpha$ contains the
composition pattern $\beta$  if there exists $i_1<i_2<\cdots<i_m$ such that
\begin{itemize}
  \item [(1)] $\alpha_{i_{s}}\leq \alpha_{i_{t}}$ if and only if $\beta_{s}\leq \beta_{t}$;
  \item [(2)] $|\alpha_{i_{s}}- \alpha_{i_{t}}|\geq |\beta_{s}-\beta_{t}|$.
\end{itemize}
We say that $\alpha$ avoids $\beta$ if $\alpha$ does not contain the pattern $\beta$.
Using the quasi-key model along with  the Kohnert diagram model of key polynomials,
Hodges and   Yong \cite{Hod} showed that $\kappa_\alpha(x)$ equals   $\Min_\alpha(x)$ if and only if $\alpha$ avoids the following five composition
patterns:
\[(0, 1, 2), (0, 0, 2, 2), (0, 0, 2, 1), (1, 0, 3, 2), (1, 0, 2, 2).\]

Our second result  gives a characterization of when   $\kappa_\alpha(x)$ reaches the upper bound.

\begin{theo}\label{main-2x}
The key polynomial $\kappa_\alpha(x)$ equals $\Max_\alpha(x)$ if and only if $\alpha$ avoids the
composition pattern $(0,2)$,
that is, there do not exist $i<j$ such that $\alpha_j-\alpha_i\ge2$.

\end{theo}

Huh,   Matherne,   M\'esz\'aros and  St.\,Dizier \cite[Conjecture 23]{Huh}   conjectured that for any
composition $\alpha\in\mathbb{Z}_{\geq 0}^n$, the normalized key polynomial $\N(\kappa_{\alpha}(x))$ is Lorentzian.
Notice that the support of $\kappa_{\alpha}(x)$ is the set of integral points in the Minkowski sum of   matroid polytopes associated to the skyline diagram $D(\alpha)$ \cite[Theorem 11]{Fin-1}.
This, along with Theorem \ref{main-2x},
enables  us to
invoke the same arguments as in the proof of   \cite[Proposition 17]{Huh}
to verify the above conjecture for compositions avoiding $(0,2)$.

\begin{coro}
If $\alpha$ avoids the
composition pattern $(0,2)$, then the key polynomial $\kappa_{\alpha}(x)$ is Lorentzian,
and hence  $\N(\kappa_{\alpha}(x))$ is Lorentzian.
\end{coro}

Note that the key polynomial $\kappa_{(0,2)}(x)$ is not Lorentzian.
This is because   $\kappa_{(0,2)}(x)$ equals the Schubert polynomial
$\S_{1423}(x)$, and the latter is not Lorentzian, as mentioned after Theorem \ref{main}.

This paper is structured  as follows. In Section \ref{Weyl}, we give an overview of
 flagged Weyl modules as well as the fact that Schubert and key polynomials are
the dual characters of   flagged Weyl modules respectively
associated to Rothe diagrams and skyline diagrams. We complete the proofs
 of Theorems \ref{main} and \ref{main-2x}
respectively in  Sections \ref{fi-pf} and \ref{fi-xx}.

\section{Dual characters of flagged Weyl modules}\label{Weyl}

Let us start with an overview of the flagged Weyl module $\mathcal{M}_D$ associated to a diagram $D$.
The module $\mathcal{M}_D$ can be  constructed  by means  of determinants \cite{Mag}.
Here we use the notation in \cite{Fin-1,Fin}.
Let $\GL(n,\mathbb{C})$ be the group of $n\times n$ invertible matrices over $\mathbb{C}$, and
let $B$ be the subgroup   consisting of the $n\times n$ upper-triangular matrices.
Let $Y$ be the $n\times n$ upper-triangular matrix whose entries are
 indeterminates $y_{ij}$ where $i\leq j$. Denote by $\mathbb{C}[Y]$
 the ring of polynomials in the variables $\{y_{ij}\}_{ i\leq j}$.
The group $\GL(n,\mathbb{C})$ acts on $\mathbb{C}[Y]$ (on the right)
as follows.  Given a matrix $g\in \GL(n,\mathbb{C})$ and a polynomial $f(Y)\in \mathbb{C}[Y]$,
\[f(Y)\cdot g=f(g^{-1} Y).\]
To a diagram
$D=(D_1,\ldots,D_n)$, the associated flagged Weyl module $\mathcal{M}_D$ is
a $B$-module  defined by
\begin{equation}\label{MM}
\mathcal{M}_D= \mathrm{Span}_\mathbb{C}\left\{\prod_{j=1}^n \det\left(Y_{D_j}^{C_j}\right)\colon C\leq D \right\},
\end{equation}
where, for two subsets $R$ and $S$ of $[n]$, $Y^R_S$ denotes the submatrix
of $Y$ with row indices in $R$ and column indices in $S$.
It should be noted  that $\prod_{j=1}^n \det\left(Y_{D_j}^{C_j}\right)\neq 0$ if and only if $C\leq D$.

Let $X=\mathrm{diag}(x_1,\ldots,x_n)$ be a  diagonal matrix, which can
be viewed as a linear transformation from $\mathcal{M}_D$ to $\mathcal{M}_D$ via the
$B$-action. The character of $\mathcal{M}_D$ is defined  as the trace of $X$:
\[\Char(\mathcal{M}_D )(x)=\tr(X\colon \mathcal{M}_D\rightarrow \mathcal{M}_D).\]
The dual character of $\mathcal{M}_D$ is the character of the dual module $\mathcal{M}_D^*$, which is given by
\begin{align*}
\Char^\ast(\mathcal{M}_D)(x)&=
\tr(X\colon \mathcal{M}_D^*\rightarrow \mathcal{M}_D^*)\\[5pt]
&=\Char(\mathcal{M}_D )(x_1^{-1},\ldots, x_n^{-1}).
\end{align*}

The Schubert and key polynomials are equal to the dual characters of
flag Weyl modules respectively corresponding to the Rothe  and skyline diagrams.
Recall that Schubert polynomials are defined based on the divided difference
operator $\partial_i$, which sends
  a polynomial $f(x)\in \mathbb{Z}[x_1,\ldots,x_n]$  to
\[\partial_i f(x)=\frac{f(x)
    -s_i f(x)}{x_i-x_{i+1}},\]
where $s_i f(x)$ is obtained from $f(x)$ by exchanging $x_i$ and $x_{i+1}$.
For the longest permutation $w_0=n \,(n-1)\cdots   1$, set
$\S_{w_0}(x)=x_1^{n-1}x_2^{n-2}\cdots x_{n-1}$.
 For $w\neq w_0$, there exists a position  $1\leq i<n$ such that $w_i<w_{i+1}$.
Let  $ws_i$ be the permutation obtained from $w$ by  interchanging
$w_i$ and $w_{i+1}$.
 Set
$\S_{w}(x)=\partial_i \S_{ws_i}(x)$.
The  above definition is independent of the choice of $i$ since the operators $\partial_i$ satisfy the
  braid relations: $\partial_i \partial_j=\partial_j \partial_i$ for $|i-j|>1$, and
$\partial_i \partial_{i+1} \partial_i= \partial_{i+1}\partial_i  \partial_{i+1}$.

As mentioned above, $\S_w(x)$ coincides with the dual character
of  $\mathcal{M}_{D(w)}$ \cite{Kra,Kra-2}:
\begin{align}\label{XY-2}
\S_w(x)= \Char^\ast(\mathcal{M}_{D(w)})(x).
\end{align}
For $C\leq D$, the effect of the action of $X$ on the polynomial $\prod_{j=1}^n \det\left(Y_{D_j}^{C_j}\right)$ is
\[\prod_{j=1}^n \det\left(Y_{D_j}^{C_j}\right)\cdot X= \prod_{j=1}^n\prod_{i\in C_j}x_i^{-1}\cdot \prod_{j=1}^n \det\left(Y_{D_j}^{C_j}\right).\]
Thus the polynomial $\prod_{j=1}^n \det\left(Y_{D_j}^{C_j}\right)$  is an
eigenvector of $X$ with eigenvalue
\[\prod_{j=1}^n\prod_{i\in C_j}x_i^{-1}.\]
Therefore, the set of monomials appearing in $\S_w(x)$ is exactly
\[\left\{x^C\colon C\le D(w)\right\}.\]
Moreover, the coefficient of a monomial $x^\alpha$
appearing in $\S_w(x)$ is equal to the dimension of the corresponding  eigenspace
\begin{align}\label{es}
\mathrm{Span}_{\mathbb{C}}\left\{\prod_{j=1}^n \det\left(Y_{D_j}^{C_j}\right)\colon  C\leq D(w), \, x^C=x^\alpha\right\}.
\end{align}
By the above observations, we obtain
the lower  and the upper bounds  for Schubert polynomials
as given in \eqref{FG-1}.

Obviously,  $\S_w(x)$ equals the lower bound $\Min_w(x)$ if and only if for each
monomial $x^\alpha$ appearing in $\S_w(x)$,   the eigenspace in \eqref{es}
has dimension one. While,
$\S_w(x)$ equals the upper bound $\Max_w(x)$ if and only if for each monomial
$x^\alpha$ appearing in $\S_w(x)$,
the    dimension of  the eigenspace in \eqref{es}
is
\[\#  \{C\colon C\leq D(w),\, x^C=x^\alpha\},\]
 that is,
the collection  of  polynomials
\[
\prod_{j=1}^n \det\left(Y_{D_j}^{C_j}\right), \quad \text{where}\ C\le D(w),
\]
are linearly independent.

Let us use  an example to illustrate \eqref{XY-2}. Consider the permutation
 $w=1432$. We have $D(1432)=(\emptyset,\{2,3\},\{2\},\emptyset)$. There are six  diagrams $C\leq D(1432)$ as listed below:
\begin{align}
&C^{(1)}=(\emptyset, \{1,2\}, \{1\},\emptyset),\ \   C^{(2)}=(\emptyset, \{1,3\}, \{1\}, \emptyset),
\ \   C^{(3)}=(\emptyset, \{2,3\}, \{1\}, \emptyset),\nonumber\\[5pt]
&C^{(4)}=(\emptyset, \{1,2\}, \{2\}, \emptyset),\ \  C^{(5)}=(\emptyset, \{1,3\}, \{2\}, \emptyset),
\ \   C^{(6)}=(\emptyset, \{2,3\}, \{2\}, \emptyset).\label{XY-1}
\end{align}
Notice  that the diagrams $C^{(3)}$ and $C^{(5)}$
give rise to the same monomial $y_{12}  y_{22}  y_{33}$.
So the module $\mathcal{M}_{D(w)}$  is spanned by the following
set of    polynomials
\[\{(y_{12}y_{23}-y_{13}y_{23})y_{12},\  (y_{12}y_{23}-y_{13}y_{23})y_{22},\ y_{12}^2y_{33},
\ y_{12} y_{22} y_{33},\  y_{22}^2y_{33} \}.\]
It is easily checked that the above five polynomials are linearly independent. So,
\begin{align*}
\Char^\ast(\mathcal{M}_{D(1432)})(x)&=x^{C^{(1)}}+x^{C^{(2)}}+x^{C^{(3)}}+x^{C^{(4)}}+x^{C^{(6)}} \\[5pt]
&=x_1^2x_2+x_1^2x_3+x_1 x_2x_3+x_1 x_2^2+x_2^2x_3,
\end{align*}
which agrees with the Schubert polynomial $\S_{1432}(x)$.

We finally  turn to key polynomials. Key polynomials  are defined using the Demazure operator $\pi_i=\partial_ix_i$.
If $\alpha$ is a partition,  set
$\kappa_\alpha(x)=x^\alpha.$
Otherwise, choose $i$ such that $\alpha_i<\alpha_{i+1}$. Let $\alpha'$ be the
composition obtained from $\alpha$ by interchanging $\alpha_i$ and $\alpha_{i+1}$. Set
$\kappa_\alpha(x)=\pi_i \kappa_{\alpha'}(x)$. The key polynomial $\kappa_\alpha(x)$ equals the dual character
of the flag Weyl module associated to the skyline diagram $D(\alpha)$:
\begin{align*}
\kappa_\alpha(x)= \Char^\ast(\mathcal{M}_{D(\alpha)})(x).
\end{align*}
In view of the  arguments for Schubert polynomials, we obtain the lower and the upper bounds
for key polynomials   given in \eqref{FG-1x}.

\section{Proof  of Theorem \ref{main}}\label{fi-pf}

In this section, we shall   prove the necessity and the sufficiency of Theorem \ref{main} in
Theorem \ref{F-1} and Theorem \ref{F-2}, respectively.

\begin{theo}\label{F-1}
If $\S_w(x)=\Max_w(x)$, then $w$ avoids the patterns 1432 and 1423.
\end{theo}

\begin{figure}[h]
\setlength{\unitlength}{0.5mm}
\begin{center}
\begin{picture}(70,100)

\put(10,30){\line(1,0){10}}\put(10,40){\line(1,0){10}}
\put(10,50){\line(1,0){10}}\put(10,60){\line(1,0){10}}
\put(10,70){\line(1,0){10}}\put(10,80){\line(1,0){10}}
\put(10,30){\line(0,1){50}}\put(20,30){\line(0,1){50}}

\qbezier[6](10,90)(15,90)(20,90)
\qbezier[6](10,80)(10,85)(10,90)
\qbezier[6](20,80)(20,85)(20,90)

\put(40,30){\line(1,0){10}}\put(40,40){\line(1,0){10}}
\put(40,30){\line(0,1){10}}\put(50,30){\line(0,1){10}}

\qbezier[6](40,50)(45,50)(50,50) \qbezier[6](40,60)(45,60)(50,60)
\qbezier[6](40,70)(45,70)(50,70) \qbezier[6](40,80)(45,80)(50,80)
\qbezier[6](40,90)(45,90)(50,90)
\qbezier[30](40,40)(40,65)(40,90)
\qbezier[30](50,40)(50,65)(50,90)

\put(-5,85){\circle*{3}}
\put(-5,85){\line(1,0){80}}\put(-5,85){\line(0,-1){90}}

\put(75,35){\circle*{3}}
\put(75,35){\line(1,0){20}}\put(75,35){\line(0,-1){40}}

\put(45,15){\circle*{3}}
\put(45,15){\line(1,0){50}}\put(45,15){\line(0,-1){20}}

\put(15,5){\circle*{3}}
\put(15,5){\line(1,0){80}}\put(15,5){\line(0,-1){10}}

\put(-17,83){\small{$i_0$}}\put(-17,33){\small{$k_0$}}
\put(-17,13){\small{$p_0$}}\put(-17,3){\small{$q_0$}}

\put(12.5,97){\small{$j_0$}}\put(42.5,97){\small{$l_0$}}
\end{picture}
\end{center}
\caption{Local configuration of the boxes in column $j_0$ and column $l_0$.}
\label{For-Pf}
\end{figure}
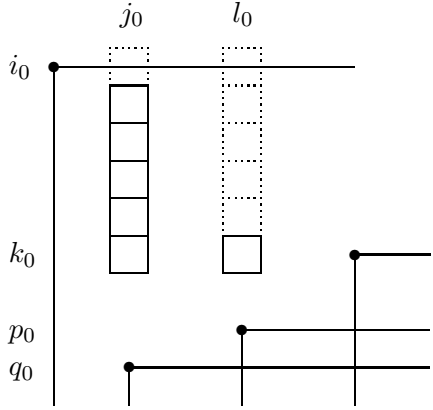

\pf We first show that if $\S_w(x)=\Max_w(x)$, then $w$ must avoid  1432.
The proof is by contradiction.
Suppose otherwise that $w$ contains a subsequence that is
order isomorphic to $1432$.
Let $i_0$ be the largest $i$ such that $w_iw_kw_pw_q$ is order isomorphic
to  1432. Once $i_0$ is determined, let $k_0$ be the smallest $k$
such that $w_{i_0}w_{k_0}w_pw_q$
is order isomorphic
to  1432. Now, let $w_{i_0}w_{k_0}w_{p_0}w_{q_0}$ be any
fixed subsequence order isomorphic
to  1432.
By the choices of $i_0$ and $k_0$, we see that $w_{q_0}<w_i<w_{p_0}$ for any $i_0<i<k_0$.
Denote  $j_0=w_{q_0}$ and $l_0=w_{p_0}$.
We have the following two  observations.
\begin{itemize}
\item[O 1.]  The box $(i_0, j_0)\notin D(w)$.
For $i_0<i\leq k_0$, the box $(i,j_0)\in D(w)$.

\item[O 2.] The box $(k_0, l_0)\in D(w)$.
For $i_0\leq i<k_0$, the box $(i,l_0)\notin D(w)$.
\end{itemize}
So the configuration of the boxes of $D(w)$ in   column $j_0$ and column $l_0$ that lie between
 row $i_0$ and row $k_0$  is as illustrated in Figure \ref{For-Pf}.

Assume that $D(w)=(D_1,\ldots, D_n)$.
For ease of description, we denote the polynomial generated by a  diagram $C\leq D(w)$
by
\begin{align}\label{PPP}
f_C(Y)=\prod_{j=1}^n \det\left(Y_{D_j}^{C_j}\right).
\end{align}
Let $t=k_0-i_0+1$. We shall construct $t$ distinct diagrams $C^{(1)}, \ldots, C^{(t)}$ such that the
corresponding polynomials $f_{C^{(1)}}(Y),\ldots, f_{C^{(t)}}(Y)$
are linearly dependent. For $1\leq m\leq t$, assume that $C^{(m)}=(C^{(m)}_1,\ldots, C^{(m)}_n)$.
The diagram   $C^{(m)}$
is defined  as follows.
\begin{itemize}
\item[(1)] For $j\notin \{j_0, l_0\}$,  let $C^{(m)}_j=D_j$.

\item[(2)] For $j=l_0$,   let
\[C^{(m)}_{l_0}=(D_{l_0}\setminus \{k_0\})\cup\{i_0+m-1\}.\]

\item[(3)] For $j=j_0$,   let
\[C^{(m)}_{j_0}=\left(D_{j_0}\cup \{i_0\}\right)\setminus \{i_0+m-1\}.\]
\end{itemize}
By the above constructions, it is easily seen that $C^{(m)}_{l_0}\leq D_{l_0}$
and $C^{(m)}_{j_0}\leq D_{j_0}$, and hence    $C^{(m)}\leq D(w)$ for each $1\leq m\leq t$.

By definition,  we have
\begin{align}\label{BN}
f_{C^{(m)}}(x)=\det\left(Y^{C^{(m)}_{j_0}}_{D_{j_0}}\right) \cdot \det\left(Y^{C^{(m)}_{l_0}}_{D_{l_0}}\right)\cdot \prod_{j\notin\{ j_0, l_0\}} \det\left(Y^{C^{(m)}_j}_{D_j}\right).
\end{align}
Now we evaluate the three factors appearing in \eqref{BN}.
For $j\notin\{ j_0, l_0\}$, since $C^{(m)}_j=D_j$, it follows  that
$\det\left(Y^{C^{(m)}_j}_{D_j}\right)$ is an upper-triangular
matrix, and thus
\begin{align}\label{lp-1}
\det\left(Y^{C^{(m)}_j}_{D_j}\right)=\prod_{j\notin\{ j_0, l_0\}} \prod_{i\in D_j} y_{ii}.
\end{align}

To calculate $\det\left(Y^{C^{(m)}_{j_0}}_{D_{j_0}}\right)$,  let
\[R^{(m)}=[i_0, k_0]\setminus \{i_0+m-1\}\]
and
\[S=D_{j_0}\setminus [i_0+1, k_0],\]
where, for two integers $a<b$, we use $[a,b]$ to denote the interval $\{a,a+1,\ldots, b\}$.
Clearly,  $C^{(m)}_{j_0}$ is the disjoint union of $R^{(m)}$ and $S$, and
$D_{j_0}$ is the disjoint union of $[i_0+1, k_0]$ and $S$.
Noticing that the submatrix of $Y^{C^{(m)}_{j_0}}_{D_{j_0}}$ obtained by restricting
both the  row and the column indices
to  $S$ is upper-triangular, we have
\begin{align}\label{lp-2}
\det\left(Y^{C^{(m)}_{j_0}}_{D_{j_0}}\right)= \det\left(Y^{R^{(m)}}_{[i_0+1, k_0]}\right)\cdot\prod_{i\in S} y_{ii}.
\end{align}

Moreover, by the choice of $C^{(m)}_{l_0}$, it is easy to see that  the matrix $Y^{C^{(m)}_{l_0}}_{D_{l_0}}$   is an upper-triangular
matrix.
So we obtain that
\begin{align}\label{lp-3}
\det\left(Y^{C^{(m)}_{l_0}}_{D_{l_0}}\right)= y_{a_m k_0}\cdot\prod_{i\in D_{l_0}\setminus \{k_0\}} y_{ii}
,
\end{align}
where $a_m=i_0+m-1$.

In view of \eqref{lp-1}, \eqref{lp-2} and \eqref{lp-3}, we find that
the polynomials $f_{C^{(1)}}(Y),\ldots, f_{C^{(t)}}(Y)$ have the following common factor
\[\prod_{j\notin \{j_0, l_0\}} \prod_{i\in D_j} y_{ii} \cdot \prod_{i\in S} y_{ii}  \cdot \prod_{i\in D_{l_0}\setminus \{k_0\}} y_{ii}.\]
Therefore, to prove that the polynomials $f_{C^{(1)}}(Y),\ldots, f_{C^{(t)}}(Y)$  are linearly dependent, it is enough to
verify that for $1\le m\le t$, the polynomials
\[g_m(Y)=y_{a_m k_0}\cdot \det\left(Y^{R^{(m)}}_{[i_0+1, k_0]}\right)\]
are linearly dependent. For simplicity, let $b=k_0$.
Without loss of generality, assume that $i_0=1$.
 Then we have $a_m=m$ and $R^{(m)}=[b]\setminus \{m\}$, and hence
\[g_m(Y)=y_{m b}\cdot \det\left(Y^{[b]\setminus \{m\}}_{[2,b]}\right),\ \ \ \text{where $1\leq m\leq b$}.\]
We claim that
\begin{align}\label{rp-1}
g_{b}(Y)=g_{b-1}(Y)-g_{b-2}(Y)+\cdots+(-1)^{b}g_1(Y).
\end{align}

To prove the claim, let us first consider
\[g_b(Y)=y_{bb}\cdot \det\left(Y^{[b-1] }_{[2,b]}\right).\]
Notice that in the  last row of the matrix $Y^{[b-1]}_{[2,b]}$,
the only nonzero entries are   $y_{(b-1)(b-1)}$ and $y_{(b-1)b}$.
Using the Laplace expansion along the last row, we have
\begin{align}\label{Mat-1}
g_b(Y)=y_{bb}\left(-y_{(b-1)(b-1)}\cdot \det\left(Y^{[b-2] }_{[2,b]\setminus \{b-1\}}\right)+
y_{(b-1)b}\cdot \det\left(Y^{[b-2] }_{[2,b-1]}\right)\right).
\end{align}
Then we consider
\[g_{b-1}(Y)=y_{(b-1)b}\cdot \det\left(Y^{[b]\setminus \{b-1\} }_{[2,b]}\right).\]
Since the last row of $Y^{[b]\setminus \{b-1\} }_{[2,b]}$ has only one nonzero element $y_{bb}$,
applying the
Laplace expansion along the last row gives
\begin{align}\label{Mat-2}
g_{b-1}(Y)=y_{(b-1)b}\cdot y_{bb}\cdot
 \det\left(Y^{[b-2] }_{[2,b-1]}\right).
\end{align}
Combining \eqref{Mat-1} and \eqref{Mat-2}, we are led to
\begin{align}\label{Mat-3}
g_b(Y)=g_{b-1}(Y)-y_{bb}\cdot y_{(b-1)(b-1)}\cdot \det\left(Y^{[b-2] }_{[2,b]\setminus \{b-1\}}\right).
\end{align}

Let us proceed to consider the summand
\[y_{bb}\cdot y_{(b-1)(b-1)}\cdot \det\left(Y^{[b-2] }_{[2,b]\setminus \{b-1\}}\right)\]
appearing in \eqref{Mat-3}. Again, applying the Laplace expansion to $Y^{[b-2] }_{[2,b]\setminus \{b-1\}}$
along the last row yields
\begin{align}\label{Mat-4}
&y_{bb}\cdot  y_{(b-1)(b-1)}\cdot \det\left(Y^{[b-2] }_{[2,b]\setminus \{b-1\}}\right)\nonumber\\[5pt]
&\quad=y_{bb}\cdot y_{(b-1)(b-1)}\left(-y_{(b-2)(b-2)} \cdot \det\left(Y^{[b-3] }_{[2,b]\setminus \{b-1, b-2\}}\right)
+y_{(b-2)b}\det\left(Y^{[b-3] }_{[2,b-2]}\right)\right).
\end{align}
On the other hand,
\begin{align}\label{Mat-5}
g_{b-2}(Y)&=y_{(b-2)b}\cdot \det\left(Y^{[b] \setminus \{b-2\}}_{[2,b]}\right)\nonumber\\[5pt]
&=y_{(b-2)b}\cdot y_{bb}\cdot \det\left(Y^{[b-1] \setminus \{b-2\}}_{[2,b-1]}\right)\nonumber\\[5pt]
&=y_{(b-2)b}\cdot y_{bb}\cdot y_{(b-1)(b-1)}\cdot \det\left(Y^{[b-3]}_{[2,b-2]}\right).
\end{align}
In view of \eqref{Mat-3}, \eqref{Mat-4} and \eqref{Mat-5}, it follows  that
\[g_b(Y)=g_{b-1}(Y)-g_{b-2}(Y)+y_{bb}\cdot y_{(b-1)(b-1)}\cdot y_{(b-2)(b-2)} \cdot \det\left(Y^{[b-3] }_{[2,b]\setminus \{b-1, b-2\}}\right).\]

Continuing the same procedure, we can arrive at the assertion in \eqref{rp-1} eventually. This
implies that the polynomials $g_m(Y)$ are linearly  dependent, and so the polynomials $f_{C^{(1)}}(Y),\ldots, f_{C^{(t)}}(Y)$ are linearly dependent.  Hence we conclude  that   $\S_w(x)\neq \Max_w(x)$.

The same arguments can be employed  to show that if $w$ contains a subsequence order-isomorphic to
1423, then   $\S_w(x)\neq \Max_w(x)$. In fact, the proof for the case 1432 does not use the relative order
of $w_{p_0}$ and $w_{q_0}$ in the subsequence  $w_{i_0}w_{k_0}w_{p_0}w_{q_0}$. This completes the proof.
\qed

We now prove  the sufficiency  of Theorem \ref{main}.

\begin{theo}\label{F-2}
If $w$ avoids the patterns 1432 and 1423, then $\S_w(x)=\Max_w(x)$.
\end{theo}

\pf Suppose that $\{C\colon C\leq D(w)\}=\{C^{(1)},\ldots, C^{(t)}\}$.
For a diagram $C$, the polynomial $f_C(Y)$ is as defined in \eqref{PPP}.
To prove  $\S_w(x)=\Max_w(x)$, it is equivalent to show   that
the polynomials $f_{C^{(1)}}(Y),\ldots, f_{C^{(t)}}(Y)$ are linearly independent.

Assume that $D(w)=(D_1,\ldots, D_n)$. For $1\leq j\leq n$, let
\[D_j'=\{i\in D_j\colon \text{there exists some $1\leq k<i$ such that $k\not\in D_j$}\}.\]
Equivalently, if we let $a_j$ be the largest integer such that $[a_j]\subseteq D_j$,
then
\[D_j'=D_j\setminus [a_j].\]
It is easy to see that for any diagram  $C=(C_1,\ldots, C_n)\leq D(w)$, we must have $[a_j]\subseteq C_j$.
This allows us to obtain the following equality
\[\det\left(Y_{D_j}^{C_j}\right)=\prod_{i=1}^{a_j}y_{ii}\cdot \det\left(Y_{D_j\setminus [a_j]}^{C_j\setminus [a_j]}\right)=\prod_{i=1}^{a_j}y_{ii}\cdot \det\left(Y_{D_j'}^{C_j\setminus [a_j]}\right),\]
and so we have
\[f_C(Y)=\prod_{j=1}^n\prod_{i=1}^{a_j}y_{ii} \cdot \prod_{j=1}^n\det\left(Y_{D_j'}^{C_j\setminus [a_j]}\right).\]
Hence, to show that the polynomials $f_{C^{(1)}}(Y),\ldots, f_{C^{(t)}}(Y)$ are linearly independent,
it suffices to show that the following polynomials  are linearly independent:
\[h_m(Y)=\prod_{j=1}^n\det\left(Y_{D_j'}^{C_j^{(m)}\setminus [a_j]}\right),\ \ \ \text{where $1\leq m\leq t$}.\]

To this end, we claim that for $1\leq j_1<j_2\leq n$, $D_{j_1}'\cap D_{j_2}'=\emptyset$.
Suppose otherwise that $D_{j_1}'\cap D_{j_2}'\neq \emptyset$.
Choose an element $i_0\in D_{j_1}\cap D_{j_2}$. By the definition of $D_j'$, there exists
$1\leq k<i_0$ that does not belong to $D_{j_1}$. This means  that the box $(k, j_1)\notin D(w)$, implying that  $w_k<j_1$.  Now we consider the subsequence $w_kw_{i_0}w_pw_q$, where
$\{w_p, w_q\}=\{j_1, j_2\}$. Since the box $(i_0,j_2)$ belongs to $D(w)$, it follows
that $w_{i_0}>j_2$. So we have $w_k<j_1<j_2<w_{i_0}$. This implies that
the subsequence $w_kw_{i_0}w_pw_q$ is order isomorphic to 1432 or 1423, leading to
a contradiction. This verifies  the claim.

By the above claim,  we see that for $1\leq t_1\neq t_2\leq t$,
 $h_{t_1}(Y)$ and $h_{t_2}(Y)$ do not contain any monomial in common,
 which  obviously   implies that the polynomials $h_m(Y)$ are linearly independent.
This completes the proof.
\qed

\section{Proof of Theorem \ref{main-2x}}\label{fi-xx}

The proof of Theorem \ref{main-2x} can be carried out in the same vein as Theorem \ref{main}.

\noindent
{\it Proof of Theorem \ref{main-2x}.}
We first prove the necessity, that is, if $\kappa_\alpha(x)=\Max_\alpha(x)$,
then  $\alpha$ avoids the
pattern $(0,2)$. Suppose otherwise that  $\alpha$ contains a
$(0,2)$ pattern. Write $\alpha=(\alpha_1,\ldots,\alpha_n)$.
Among the $(0,2)$ patterns of $\alpha$, choose the largest index
$i_0$ such that $(\alpha_{i_0}, \alpha_k)$ is a $(0,2)$ pattern for some $k>i_0$.
Once $i_0$ is fixed, locate the smallest index $k_0$
such that $(\alpha_{i_0}, \alpha_{k_0})$ is a $(0,2)$ pattern.
Since $(\alpha_{i_0}, \alpha_{k_0})$ is a $(0,2)$ pattern,
by definition we have
\begin{equation}\label{FN-1}
\alpha_{k_0}-\alpha_{i_0}\geq 2.
\end{equation}
By the choices of $i_0$ and $k_0$, it is easy to check that
for $i_0<i<k_0$
\begin{equation}\label{FN-2}
\alpha_i-\alpha_{i_0}=1.
\end{equation}

\begin{figure}[h]
\begin{center}
\begin{tikzpicture}

\def\rectanglepath{-- +(4mm,0mm) -- +(4mm,4mm) -- +(0mm,4mm) -- cycle}


\draw [step=4mm] (0mm,20mm) grid (12mm,32mm);
\draw(0mm,20mm)--(12mm,20mm);
\draw (12mm,20mm) \rectanglepath;
\draw (12mm,24mm) \rectanglepath;

\node at (6mm,17mm) {$\vdots$};

\draw [step=4mm] (0mm,4mm) grid (16mm,12mm);
\draw (20mm,4mm) \rectanglepath;\draw (16mm,4mm) \rectanglepath;
\draw (24mm,4mm) \rectanglepath;

\node at (-4mm,6mm) {$k_0$};\node at (-4mm,30mm) {$i_0$};
\node at (15mm,36mm) {$j_0$};\node at (26mm,36mm) {$l_0$};

\end{tikzpicture}
\end{center}
\vspace{-4mm}
\caption{An illustration of the  boxes of $D(\alpha)$ in row $i_0$ and $k_0$.}
\label{MMM}
\end{figure}
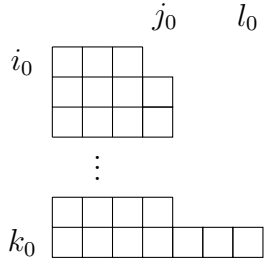

 By \eqref{FN-1}
and \eqref{FN-2}, the configuration of the boxes of $D(\alpha)$ lying  between
 row $i_0$ and row $k_0$  is depicted in Figure \ref{MMM}.
Let $j_0=\alpha_{i_0}+1$ and $l_0=\alpha_{k_0}$.
Clearly,
the configuration of the boxes of $D(\alpha)$ in   column $j_0$ and column $l_0$ that lie between
 row $i_0$ and row $k_0$   is completely that same as that in Figure \ref{For-Pf}.
 Therefore, using the same arguments as in the proof of  Theorem \ref{F-1},
 we  obtain that the assumption that $\alpha$ contains a
$(0,2)$ pattern is false. This verifies the necessity.

It remains to prove the sufficiency. The analysis  is  similar to that in
the proof of Theorem \ref{F-2}, and is sketched below. Assume that
$\alpha$ avoids the composition pattern $(0,2)$.
Write $D(\alpha)=(D_1,\ldots, D_n)$. For $1\leq j\leq n$,
let $a_j$ be the largest integer such that $[a_j]\subseteq D_j$.
Define
\[D_j'=D_j\setminus [a_j].\]
For $1\leq j_1<j_2\leq n$, we claim that $D_{j_1}'\cap D_{j_2}'=\emptyset$.
Suppose to the contrary that $D_{j_1}'\cap D_{j_2}'\neq\emptyset$.
Assume that  $t\in D_{j_1}'\cap D_{j_2}'$. Denote $s=a_{j_1}+1$.
By the definition of $a_j$, we see that $s$ is smaller than any
integer appearing in $D_{j_1}'$. Since $t\in D_{j_1}'$, we have $s<t$.
 Consider the parts $\alpha_s$ and $\alpha_t$.
 Since $D(\alpha)$ is a skyline diagram, it is clear that $D_{j_1}\supseteq D_{j_2}$,
which, along with the fact that $s\not\in D_{j_1}$, implies    $s\not\in D_{j_2}$. Combining the
assumption that  $t\in D_{j_1}\cap D_{j_2}$, we see that $\alpha_s\leq \alpha_t-2$,
and so  $(\alpha_s,\alpha_t)$ forms a $(0,2)$ patten,
leading to a
contradiction. This verifies   the claim.

Based on the above claim, we can now use the same arguments as in the
proof of Theorem \ref{F-2} to conclude the  sufficiency.
This finishes the proof.
\qed

\vspace{.2cm} \noindent{\bf Acknowledgments.}
This work was supported by the National Natural
Science Foundation of China (Grant No. 11971250)
and Sichuan Science and Technology Program (Grant No. 2020YJ0006).


\begin{thebibliography}{99}

\bibitem{Ard}
F. Ardila,
The Catalan matroid,
 J. Combin. Theory Ser. A 104 (2003),  49--62.

\bibitem{Ass}
 S. Assaf,
Nonsymmetric Macdonald polynomials and a refinement of
Kostka--Foulkes polynomials,
Trans. Amer.  Math. Soc.  370 (2018), 8777--8796.


\bibitem{Ber}
N. Bergeron and S. Billey, RC-graphs and Schubert polynomials,
Experiment. Math. 2 (1993), 257--269.

\bibitem{Bil}
S. Billey, W. Jockusch and R.P. Stanley, Some combinatorial properties of Schubert
polynomials, J. Algebraic Combin. 2 (1993), 345--374.


\bibitem{Bon}
J. Bonin, A. de Mier and M. Noy,
Lattice path matroids: enumerative aspects and Tutte polynomials,
J. Combin. Theory Ser. A 104 (2003),  63--94.

\bibitem{Bra}
P. Br\"and\'en and J. Huh,
Lorentzian polynomials,
Ann. Math., to appear.


\bibitem{Dem-1}
M. Demazure,
D\'esingularisation des vari\'et\'es de Schubert g\'en\'eralis\'ees,
Ann. Sci. \'Ecole  Norm. Sup.  7 (1974), 53--88.

\bibitem{Dem-2}
M. Demazure,
Une nouvelle formule des caract\'eres,
Bull. Sci. Math. 98 (1974), 163--172.

\bibitem{Fan}
N.J.Y. Fan and P.L. Guo,
Vertices of Schubitopes,
J. Combin. Theory Ser. A, to appear.

\bibitem{Fin-1}
A. Fink, K. M\'esz\'aros and A. St.\,Dizier,
Schubert polynomials as integer point transforms of generalized permutahedra,
Adv. Math.  332 (2018), 465--475.



\bibitem{Fin}
A. Fink, K. M\'esz\'aros and A. St.\,Dizier,
Zero-one Schubert polynomials,
Math. Z., to appear.


\bibitem{Gir}
S. Gire, Arbres, permutations \'a motifs exclus et cartes planaire: quelques probl\'emes algorithmiques et
combinatoires, Ph.D. Thesis, University of Bordeaux, 1993.


\bibitem{Hod}
R. Hodges and A. Yong,
Multiplicity-free key polynomials,
 arXiv:2007.09229v1.


\bibitem{Huh}
J. Huh, J.P. Matherne, K. M\'esz\'aros and A. St.\,Dizier,
Logarithmic concavity of Schur and related polynomials,
arXiv:1906.09633v3.


\bibitem{Ion}
B. Ion,  Nonsymmetric Macdonald polynomials and Demazure characters,
 Duke Math. J. 116 (2003),  299--318.

\bibitem{Knu}
A. Knutson and E. Miller, Gr\"obner geometry of Schubert polynomials, Ann. Math.
161 (2005), 1245--1318.


\bibitem{Kra}
W. Kra\'skiewicz and P. Pragacz, Foncteurs de Schubert,
C. R. Acad. Sci. Paris  S\'er. I Math. 304 (1987), 209--211.

\bibitem{Kra-2}
W. Kra\'skiewicz and P. Pragacz,
Schubert functors and Schubert polynomials,
European  J. Combin.  25 (2004), 1327--1344.

\bibitem{Kremer}
D. Kremer,
Permutations with forbidden subsequences
and a generalized Schr\"{o}der number, Discrete Math. 218 (2000), 121--130.



\bibitem{Lam}
 T. Lam, S. Lee and M. Shimozono, Back stable Schubert calculus,
arXiv:1806.11233v1.

\bibitem{Las-1}
A. Lascoux and M.-P. Sch\"utzenberger, Keys \& standard bases, Invariant Theory and Tableaux (Minneapolis, MN, 1988), 125--144, IMA Vol. Math. Appl., 19, Springer, New York, 1990.


\bibitem{Las}
A. Lascoux and M.-P. Sch\"utzenberger, Polyn$ \hat{\mathrm o}$mes de Schubert, C. R. Acad. Sci.
Paris  S\'er. I Math. 294 (1982), 447--450.

\bibitem{Las-2}
A. Lascoux and M.-P. Sch\"utzenberger, Tableaux and non-commuative Schubert
polynomials, Func. Anal. Appl. 23 (1989), 63--64.



\bibitem{Mag}
 P. Magyar, Schubert polynomials and Bott-Samelson varieties,
 Comment. Math. Helv.  73 (1998), 603--636.




\bibitem{Rei}
V. Reiner and M. Shimozono,
Key polynomials and a flagged
Littlewood-Richardson rule,
J. Combin. Theory  Ser. A 70 (1995), 107--143.

\bibitem{Win}
R. Winkel, Diagram rules for the generation of Schubert polynomials, J. Combin.
Theory Ser. A. 86 (1999), 14--48.


\end{thebibliography}
\end{document}